\documentclass[11pt,reqno,oneside]{amsart}

\usepackage{graphicx}    
\usepackage{amsfonts}
\usepackage{verbatim}
\usepackage{amssymb}
\usepackage{amsmath}
\usepackage{hyperref}
\usepackage[numbers,sort&compress]{natbib}

\usepackage[it,small]{caption}  

\newcommand{\R}{{\mathbb R}}

\newcommand{\E}{\mathcal{E}}

\newcommand{\half}{{\textstyle \frac{1}{2}}}
\newcommand{\fourth}{{\textstyle \frac{1}{4}}}

\renewcommand{\vec}[1]{\mathbf{#1}}

\newtheorem{theorem}{Theorem}[section]

\numberwithin{equation}{section}

\allowdisplaybreaks[4]

\begin{document}

\title[Goal-Oriented Atomistic-Continuum Adaptivity for the Quasicontinuum Approximation]
{Goal-oriented Atomistic-Continuum Adaptivity for the Quasicontinuum Approximation}
\author{Marcel Arndt}
\author{Mitchell Luskin}
\begin{abstract}
  We give a goal-oriented {\em a posteriori} error estimator for the atomistic-continuum modeling
  error in the quasicontinuum method, and we use this estimator to design an adaptive algorithm to
  compute a quantity of interest to a given tolerance by using a nearly minimal number of atomistic
  degrees of freedom.  We present computational results that demonstrate the effectiveness of our
  algorithm for a periodic array of dislocations described by a Frenkel-Kontorova type model.
\end{abstract}



\thanks{This work was supported in part by DMS-0304326 and by the Minnesota
  Supercomputing Institute. This work is also based on work supported by the
  Department of Energy under Award Number DE-FG02-05ER25706.}

\date{\today}

\maketitle


\section{Introduction}  \label{SecIntro}

Multiscale methods offer the potential to
solve complex problems by utilizing a fine-scale model
 only in regions that require increased accuracy.  However,
 it is usually not known {\it a priori} which region requires
 the increased accuracy, and an adaptive algorithm is needed
 to compute a given quantity of interest to a given tolerance
 with nearly optimal computational efficiency.

The quasicontinuum (QC) method
\cite{TadmorMillerPhillipsOrtiz:1999, TadmorOrtizPhillips:1996,
TadmorPhilipsOrtiz:1996,DobsonLuskin:2006} is a multiscale computational method
for crystals that retains sufficient
accuracy of the atomistic model by utilizing a strain energy
density obtained from the atomistic model by the Cauchy-Born rule
in regions where the deformation is nearly uniform.
The atomistic model is needed to accurately model the
material behavior in regions of highly non-uniform deformations
around defects such as dislocations and cracks.

The approximation error within the quasicontinuum method can be decomposed into the modeling error
which occurs when replacing the atomistic model by a continuum model, and the coarsening error which
arises from reducing the number of degrees of freedom within the continuum region. This paper purely
focuses on the estimation of the modeling error. The optimal strategy to determine the mesh size in
the continuum region will be studied in a forthcoming paper.

The development of goal-oriented error estimators for mesh coarsening in the
quasicontinuum method has been given in
\cite{OdenPrudhommeRomkesBauman:2006, OdenPrudhommeBauman:2006}, and
goal-oriented error estimators for atomistic-continuum modeling have
recently been given in \cite{ArndtLuskin:2007a}.
In all these
works, the error is measured in
terms of a user-definable quantity of interest instead of a global norm.  Goal-oriented error
estimation in general is based on duality techniques and has already been used for a variety of
applications such as mesh refinement in finite element methods \cite{AinsworthOden:2000,
BangerthRannacher:2003} and control of the modeling error in homogenization \cite{OdenVemaganti:2000}.

In \cite{ArndtLuskin:2007a}, an {\em a posteriori} error estimator
for the modeling error in the quasicontinuum method has been
developed, analyzed, and applied to a one-dimensional Frenkel-Kontorova model
of crystallographic defects \cite{Marder:2000}. In this paper, we
summarize this approach and adapt it to a different setting.
Instead of clamped boundary conditions, we use periodic boundary
conditions here which are physically more realistic and allow for
more succinct formulas.  Moreover, an asymmetric quantity of
interest is used here rather than the symmetric one
studied in \cite{ArndtLuskin:2007a} to give further insight into
the behavior of the error estimator.

The one-dimensional periodic Frenkel-Kontorova model chosen here allows for an easy study of the error
estimator and keeps the formulas short, but still exhibits enough of the features of
multidimensional problems for a realistic
study of the error estimator.  In addition to the nearest-neighbor harmonic interactions between the
atoms in the classical Frenkel-Kontorova model, we add next-nearest-neighbor harmonic interactions
to obtain a non-trivial quasicontinuum approximation.

We describe the atomistic model and
its quasicontinuum approximation in Section~\ref{SecModel}, and we formulate the error estimator in
Section~\ref{SecErrorEst}. We then develop
in Section~\ref{SecErrorEst} an algorithm which employs
the error estimator for an adaptive strategy, and we conclude by
exhibiting and interpreting  the results of our numerical
experiments.

\section{Atomistic Model and Quasicontinuum Model}  \label{SecModel}

\begin{figure}
\includegraphics[width=\textwidth]{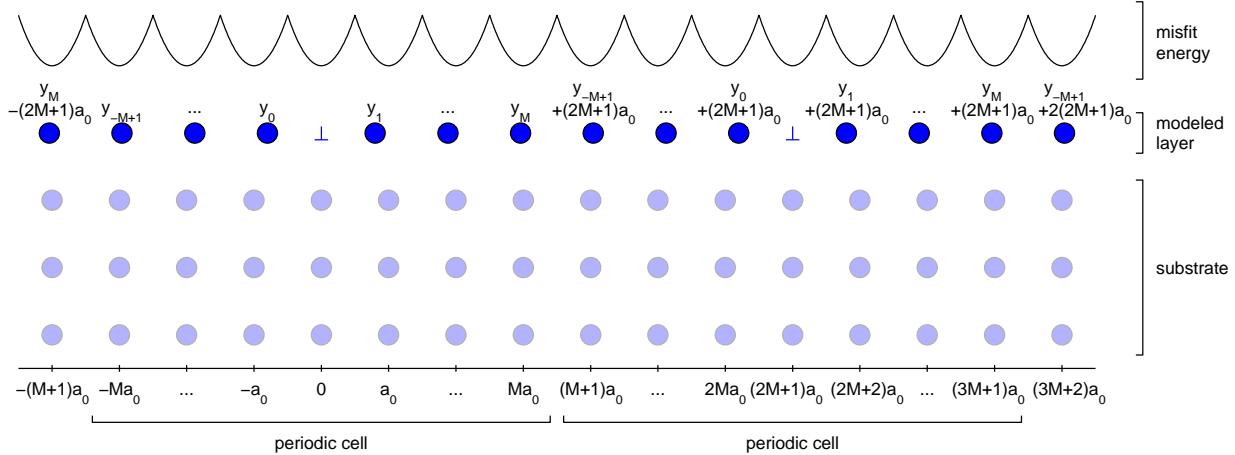}
\caption{Frenkel-Kontorova model. The wells depict the misfit
energy caused by the substrate.}
\label{FigFKModel}
\end{figure}

As an application for the method of error estimation described in
this paper, we study a periodic array of dislocations in a single
crystal. We employ a Frenkel-Kontorova type model
\cite{Marder:2000} to give a simplified one-dimensional
description of these typically two-dimensional or
three-dimensional phenomena. To model the elastic interactions, we
consider $2M$ atoms in a periodic chain that interact by Hookean
nearest-neighbor and next-nearest-neighbor springs.  The
dislocation is modeled by the interaction with a substrate which
gives rise to a misfit energy, see Figure~\ref{FigFKModel}.

The vector $\vec y = (y_{-M+1},\dots, y_M)\in\R^{2M}$
describes the positions of $2M$ atoms
which generate the positions of
an infinite chain of atoms by the relation
\begin{equation}\label{period}
y_{i+2M}=y_i+(2M+1)a_0\quad\text{for }i=-\infty,\dots,\infty,
\end{equation}
where $a_0$ denotes the lattice constant.  The relation
\eqref{period} gives an average strain of $(2M+1)/2M$ due to
a periodic array of dislocations that stretch the chain by
one lattice constant every $2M$ atoms.

The total energy $\E^a$ for this atomistic system reads
\begin{align} \label{EqFK}
\E^{a}(\vec y) = \E^{a,e}(\vec y) + \E^{a,m}(\vec y)
\end{align}
with elastic energy given by (recall \eqref{period})
\begin{equation}\label{EqFKElastic}
\begin{split}
\E^{a,e}(\vec y)&=
\half k_1  \sum_{i=-M+1}^{M} (y_{i+1} - y_i - a_0)^2
+ \half k_2  \sum_{i=-M+2}^{M+1} (y_{i+1} - y_{i-1} - 2a_0)^2\\
& = \half k_1 \Bigg[ \sum_{i=-M+1}^{M-1} (y_{i+1} - y_i - a_0)^2
                    + (y_{-M+1} - y_M + 2M a_0)^2 \Bigg] \\
& \quad + \half k_2 \Bigg[ \sum_{i=-M+2}^{M-1} (y_{i+1} - y_{i-1} - 2a_0)^2
                    + (y_{-M+1} - y_{M-1} + (2M-1)a_0)^2 \\
& \hspace{63.5mm}     + (y_{-M+2} - y_{M} + (2M-1)a_0)^2 \Bigg]
\end{split}
\end{equation}
and misfit energy given by
\begin{align} \label{EqFKMisfit1}
\E^{a,m}(\vec y)
& = \half k_0 \sum_{i=-M+1}^M \left( y_i - a_0 \left\lfloor
   \frac{y_i}{a_0} + \frac{1}{2} \right\rfloor \right)^2.
\end{align}
Here $k_0$, $k_1$ and $k_2$ denote the elastic constants.
In the misfit energy, $\lfloor x \rfloor$ denotes the largest
integer smaller than or equal to
$x.$
We can obtain the following more symmetric form of the elastic energy \eqref{EqFKElastic}
by realizing that the forces of constraint corresponding to the strain \eqref{period}
move the equilibrium spacing of the chain to $\frac{2M+1}{2M}a_0:$
\begin{equation}\label{symm}
\begin{split}
\E^{a,e}(\vec y)&=
\half k_1 \left[ \sum_{i=-M+1}^{M}
\Big( y_{i+1} - y_i - {\textstyle \frac{2M+1}{2M}}
\big(\left(\left(i+1\right) \text{ mod } 2M\right)-i\big) a_0 \Big)^2
+ {\textstyle \frac{1}{2M}}a_0^2 \right]\\
\qquad&+\half k_2 \left[ \sum_{i=-M+1}^{M} \Big( y_{i+2} - y_i - {\textstyle \frac{2M+1}{2M}}
 \big(\left(\left(i+2\right) \text{ mod } 2M\right)-i\big) a_0 \Big)^2
+   {\textstyle \frac{2}{M}}a_0^2 \right]
\end{split}
\end{equation}
where all atom indices are understood modulo $2M$, and $i \text{ mod } 2M$
denotes $i$ modulo $2M$ within the interval $-M+1,\ldots,M$. In the following,
we neglect the constant terms since they do not have any effect when
finding energy minimizers later.

We consider a single vacancy between the atoms $y_0$ and $y_1$. If we assume
that the $M$ leftmost atoms $y_i$ for $-M+1 \le i \le 0$ are situated in the
interval $\left( \left( i-\frac{3}{2} \right) a_0, \left( i-\frac{1}{2} \right)
  a_0 \right)$ and that the $M$ rightmost atoms $y_i$ for $1 \le i \le M$ are
situated in the interval $\left( \left( i-\frac{1}{2} \right) a_0, \left(
    i+\frac{1}{2} \right) a_0 \right)$, then the misfit energy can be
    expressed more simply as
\begin{align} \label{EqFKMisfit2}
\E^{a,m}(\vec y)
& = \half k_0 \sum_{i=-M+1}^M \left( y_i - b_i \right)^2
\qquad \text{with} \qquad
b_i = \begin{cases}
  (i-1)a_0 & \text{for } i \le 0, \\
  i    a_0 & \text{for } i \ge 1.
\end{cases}
\end{align}

To reduce the work
in computing \eqref{EqFK}, we employ the
quasicontinuum method \cite{TadmorMillerPhillipsOrtiz:1999,
  TadmorOrtizPhillips:1996, TadmorPhilipsOrtiz:1996} which has been successfully
used for many applications. The quasicontinuum method
consists of two steps: the passage to a continuum
energy within the continuum region of the chain, and a subsequent coarsening
in the continuum region to reduce the number of degrees of freedom.

In the first step, we replace the atomistic energy of all atoms from the
continuum region by a continuum energy. To this end, we split the total energy
into atom-wise contributions:
\begin{align}
\E^a(\vec y) = \sum_{i=-M+1}^M \E^{a,i}(\vec y)
\end{align}
with
\begin{equation}
\begin{split}
\E^{a,i}(\vec y) = \quad \fourth k_1 \Big[ \quad
  & \left( y_i - y_{i-1} - {\textstyle \frac{2M+1}{2M}} \big(i-\left(\left(i-1\right) \text{ mod } 2M\right)\big) a_0 \right)^2
    \\
+ & \left( y_{i+1} - y_i - {\textstyle \frac{2M+1}{2M}} \big(\left(\left(i+1\right) \text{ mod } 2M\right)-i\big) a_0 \right)^2 \Big] \\
+ \fourth k_2 \Big[ \quad
  & \left( y_i - y_{i-2} - {\textstyle \frac{2M+1}{2M}} \big(i-\left(\left(i-2\right) \text{ mod } 2M\right)\big) a_0 \right)^2
    \\
+ & \left( y_{i+2} - y_i - {\textstyle \frac{2M+1}{2M}} \big(\left(\left(i+2\right) \text{ mod } 2M\right)-i\big) a_0 \right)^2 \Big] \\
+ \half k_0 ( y_i & - b_i )^2.
\end{split}
\end{equation}

The corresponding continuum
energy can be derived following
\cite{ArndtLuskin:2007a} to be
\begin{equation}
\begin{split}
\E^{c,i}(\vec y) = \quad \left( \fourth k_1 + k_2 \right) \Big[ \quad
  & \left( y_i - y_{i-1} - {\textstyle \frac{2M+1}{2M}} \big(i-\left(\left(i-1\right) \text{ mod } 2M\right)\big) a_0 \right)^2
    \\
+ & \left( y_{i+1} - y_i - {\textstyle \frac{2M+1}{2M}} \big(\left(\left(i+1\right) \text{ mod } 2M\right)-i\big) a_0 \right)^2 \Big] \\
+ \half k_0 ( y_i & - b_i )^2.
\end{split}
\end{equation}
 If
\begin{align}
\delta_i^a = \begin{cases}
  1 & \text{if atom $i$ is modeled atomistically,} \\
  0 & \text{if atom $i$ is modeled as continuum,}
\end{cases}
\qquad \text{and} \qquad \delta_i^c = 1-\delta_i^a,
\end{align}
then
\begin{align}
  \E^{ac}(\vec y) = \sum_{i=-M+1}^M \left[
    \delta^a_i \E^{a,i}(\vec y) + \delta^c_i \E^{c,i}(\vec y) \right]
\end{align}
denotes the mixed atomistic-continuum energy.

In the second step of the quasicontinuum approximation, the chain is coarsened
in the continuum region by choosing representative atoms, more briefly called
repatoms.  The chain is then fully modeled in terms of the repatoms. The missing
atoms are implicitly reconstructed by linear interpolation according to the
Cauchy-Born hypothesis.  The lengthy expression of the resulting quasicontinuum
energy
\begin{align}
  \E^{qc}(\vec y)
\end{align}
is not needed in this paper since we focus on the estimation of the modeling
error. Hence we refer to \cite{ArndtLuskin:2007a} for
the formula and its derivation.  The error arising from coarsening will be
studied in a forthcoming paper.

For the subsequent argumentation, it is useful to rewrite the energies in matrix
notation. We have
\begin{equation}
\begin{split}
\E^a(\vec y)   & = \half (\vec y-\vec a^a)^T D^T E^a D (\vec y-\vec a^a)
 + \half (\vec y-\vec b^a)^T K^a (\vec y-\vec b^a), \\
\E^{ac}(\vec y) & = \half (\vec y-\vec a^a)^T D^T E^{ac} D (\vec y-\vec a^a)
 + \half (\vec y-\vec b^a)^T K^a (\vec y-\vec b^a),
\end{split}
\end{equation}
where the $2M\times 2M$ matrices are given by
\begin{equation}
\begin{array}{ll}
  D_{i,i} = -1, &
  D_{i,i+1} = 1, \\
  (E^a)_{i,i} = k_1 + 2 k_2, \qquad &
  (E^a)_{i,i+1} = (E^a)_{i+1,i} = k_2, \\
  \multicolumn{2}{l}{(E^{ac})_{i,i}
  = \half k_1 \left( \delta^a_i + \delta^a_{i+1} \right)
  + \half k_2 \left( \delta^a_{i-1} + \delta^a_{i} + \delta^a_{i+1} + \delta^a_{i+2} \right)
  + \left( \half k_1 + 2 k_2 \right) \left( \delta^c_i + \delta^c_{i+1} \right),} \\
  \multicolumn{2}{l}{(E^{ac})_{i,i+1} = (E^{ac})_{i+1,i}
  = \half k_2 \left( \delta^a_i + \delta^a_{i+2} \right),} \\
  (K^a)_{i,i} = k_0,
\end{array}
\end{equation}
with $i=-M+1, \ldots, M$ and all indices to be understood modulo $2M$ as before.
The vectors $\vec a^a\in\R^{2M}$ and  $\vec b^a\in\R^{2M}$ are defined as
\begin{equation}
\begin{split}
  \vec a^a & =
  \begin{bmatrix} (-M+1) {\textstyle \frac{2M+1}{2M}} a_0 &
                  (-M+2) {\textstyle \frac{2M+1}{2M}} a_0 &
                  \cdots &
                  ( M-1) {\textstyle \frac{2M+1}{2M}} a_0 &
                    M    {\textstyle \frac{2M+1}{2M}} a_0 \end{bmatrix}^T, \\
  \vec b^a & =
  \begin{bmatrix} b_{-M+1} & b_{-M+2} & \cdots & b_{M-1} & b_M \end{bmatrix}^T.
\end{split}
\end{equation}

We require that the elastic moduli satisfy
$k_1+2k_2>2|k_2|$ to ensure that $E^a$ is positive definite
and that the misfit modulus $k_0>0$ to ensure that
$K^a$ is positive definite.

We are interested in finding energy minimizing configurations $\vec y^a$, $\vec
y^{ac}$, and $\vec y^{qc}$ of $\E^a$, $\E^{ac}$, and $\E^{qc}$, respectively.
The minimizers $\vec y^a$ and $\vec y^{ac}$ satisfy the linear equations
\begin{equation} \label{EqPrimal}
\begin{split}
  M^a    \vec y^a   & = \vec f^a,   \\
  M^{ac} \vec y^{ac} & = \vec f^{ac}, 
\end{split}
\end{equation}
where
\begin{equation} \label{EqMfDef}
\begin{aligned}
  M^a         & := D^T E^a D + K^a, &   \hspace{20mm}
  \vec f^a    & := D^T E^a D \vec a^a + K^a \vec b^a, \\
  M^{ac}      & := D^T E^{ac} D + K^a, &
  \vec f^{ac} & := D^T E^{ac} D \vec a^a + K^a \vec b^a. 
\end{aligned}
\end{equation}
We refer to \eqref{EqPrimal} as the primal problems.
Note that the minimizers are uniquely determined due to the convexity of
the energy.

\section{Error Estimation}  \label{SecErrorEst}

In the preceding section, we described how the quasicontinuum
method gives an approximation $\vec y^{ac}$ of the atomistic
solution $\vec y^a$ by passing from the fully atomistic model to a
mixed atomistic-continuum formulation, and then we briefly
mentioned how a further approximation, $\vec y^{qc}$, can be
obtained by coarsening in the continuum region.

Instead of measuring the error in some global norm, we measure the error of a
quantity of interest denoted by $Q(\vec y)$ for some function $Q: \R^{2M}\to\R$.
We assume that $Q$ is linear and thus has a representation
\begin{align}
  Q(\vec y) = \vec q^T \vec y
\end{align}
for some vector $\vec q \in \R^{2M}$.  We then have the splitting
\begin{align}
  | Q(\vec y^a) - Q( \vec y^{qc}) |
  =   |Q(\vec y^a - \vec y^{ac})  +  Q(\vec y^{ac} - \vec y^{qc})|
  \le |Q(\vec y^a - \vec y^{ac})| + |Q(\vec y^{ac} - \vec y^{qc})|
\end{align}
of the total error into the modeling error, $|Q(\vec y^a - \vec y^{ac})|$, and the
coarsening error, $|Q(\vec y^{ac} - \vec y^{qc})|$, everything in terms of the
quantity of interest.  In this paper, we restrict ourselves to the estimation of
the modeling error.  The coarsening error will be analyzed in a forthcoming
paper.

An important tool for the error estimation in terms of a quantity of interest
are the dual problems for the influence or generalized Green's functions
$\vec g^a$ and $\vec g^{ac}$ given by
\begin{equation} \label{EqDual}
\begin{split}
  M^{a}  \vec g^a   & = \vec q, \\
  M^{ac} \vec g^{ac} & = \vec q.
\end{split}
\end{equation}
The matrices $M^a$ and $M^{ac}$ are symmetric since they stem from an
energy, and we thus do not need to use their transpose for the dual problems.

We denote the errors and the residuals, both for the deformation $\vec y^a$
and the influence function $\vec g^a,$ by
\begin{equation} \label{EqErrorPrimalDual}
\begin{aligned}
  \vec e          & := \vec y^a - \vec y^{ac}, & \hspace{20mm}
  R(\vec y)       & := M^a \left( \vec y^a - \vec y \right)
                     = \vec f^a - M^a \vec y, \\
  \hat{\vec e}    & := \vec g^a - \vec g^{ac}, &
  \hat{R}(\vec g) & := M^a \left( \vec g^a - \vec g \right)
                     = \vec q - M^a \vec g.
\end{aligned}
\end{equation}
Then we have the basic identity for the error of the quantity of
interest
\begin{equation}\label{EqBasicDual}
\begin{split}
  Q(\vec y^a) - Q(\vec y^{ac})
  & = \vec q^T \vec e = \vec g^{aT} M^a \vec e
    = (\vec g^{ac} + \hat{\vec e})^T M^a \vec e \\
  & = \vec g^{acT} R(\vec y^{ac}) + \hat{\vec e}^T M^a \vec e.
\end{split}
\end{equation}

The quantities $\vec y^{ac}$ and $\vec g^{ac}$ are considered to be computable since
the continuum degrees of freedom give local interactions, whereas $\vec y^a$
and $\vec g^a$ are not considered to be computable since they require
a full atomistic computation.  Thus the first term $\vec g^{acT} R(\vec y^{ac})$ is
easily computable, and the challenge is to estimate $\hat{\vec e}^T M^a \vec e$.
Let us note that in applications to mesh refinement for linear finite elements,
the residual term vanishes due to Galerkin orthogonality, whereas in other
applications it can be dominant over the second term.

We utilize two error estimators derived in
\cite{ArndtLuskin:2007a} and briefly summarized here.
Our first error estimator is based on the generalized parallelogram identity
\begin{equation}\label{EqParallelogram}
  \hat{\vec e}^T M^a \vec e
= \textstyle \frac{1}{4} \| \sigma\vec e + \sigma^{-1}\hat{\vec e} \|_{M^a}^2
               - \frac{1}{4} \| \sigma\vec e - \sigma^{-1}\hat{\vec e} \|_{M^a}^2
\end{equation}
for all $\sigma\ne0$, where the $M^a$-norm of some vector $\vec z$ is defined by
$\|\vec z\|_{M^a} := (\vec z^T M^a \vec z)^{1/2}$.  We define the computable
bounds
\begin{equation}\label{eta}
  \eta_\text{low}^\pm \le \| \sigma\vec e \pm \sigma^{-1}\hat{\vec e} \|_{M^a} \le \eta_\text{upp}^\pm
\end{equation}
by
\begin{equation}
\begin{aligned}
  \eta_\text{upp}^\pm
  & := \big\| P D \big[ \sigma ( \vec y^{ac} - \vec a^a )
       \pm \sigma^{-1} \vec g^{ac} \big] \big\|_{E^a}, \\
  \eta_\text{low}^\pm
  & := \frac{\left| (\vec y^{ac} + \theta^\pm \vec g^{ac})^T \vec r^\pm \right|}
       {\| \vec y^{ac} + \theta^\pm \vec g^{ac} \|_{M^a}}
\end{aligned}
\end{equation}
where
\begin{equation}
\begin{aligned}
  P          & := I - (E^a)^{-1} E^{ac}, \\
  \vec r^\pm & := \sigma R(\vec y^{ac}) \pm \sigma^{-1} \hat R(\vec g^{ac}).
\end{aligned}
\end{equation}
Optimization of the bounds with respect to $\sigma$ and $\theta$ leads
in \cite{ArndtLuskin:2007a} to the following choice of the
parameters:
\begin{equation}
\begin{split}
  \sigma & := \sqrt{\frac{\| P D \vec g^{ac} \|_{E^a}}
       {\| P D ( \vec y^{ac} - \vec a^a ) \|_{E^a}}}, \\
  \theta^\pm & :=
    \frac{ \vec r^{\pm T} \vec y^{ac} \; \vec g^{acT} M^a \vec y^{ac} -
           \vec r^{\pm T} \vec g^{ac} \; \|\vec y^{ac}\|_{M^a}^2}
         { \vec r^{\pm T} \vec g^{ac} \; \vec g^{acT} M^a \vec y^{ac} -
           \vec r^{\pm T} \vec y^{ac} \; \|\vec g^{ac}\|_{M^a}^2}.
\end{split}
\end{equation}

\begin{theorem} \label{Theorem1}
  We have that
  \begin{align}
    \left| Q(\vec y^a) - Q(\vec y^{ac}) \right| \le \eta_1,
  \end{align}
  where the computable error estimator $\eta_1$ is defined as
  \begin{align}  \label{EqEta1Def}
    \textstyle \eta_1 :=
    \max\left( \left| \vec g^{acT} R(\vec y^{ac}) + \frac{1}{4} (\eta_\text{low}^+)^2
                                             - \frac{1}{4} (\eta_\text{upp}^-)^2 \right|,
      \left| \vec g^{acT} R(\vec y^{ac}) + \frac{1}{4} (\eta_\text{upp}^+)^2
                                             - \frac{1}{4} (\eta_\text{low}^-)^2 \right| \right).
  \end{align}

\end{theorem}

We also developed the following weaker estimator in \cite{ArndtLuskin:2007a}
using the Cauchy-Schwarz inequality
in place of the parallelogram identity in \eqref{EqParallelogram}.
We note that this estimator can be decomposed among the degrees of freedom
and can thus be utilized in adaptive atomistic-continuum modeling decisions.
\begin{theorem}
  We have that
  \begin{align}
    \left| Q(\vec y^a) - Q(\vec y^{ac}) \right| \le \eta_2
    \le \sum_{i=-M+1}^{M} \eta_{2,i}^{at} + \sum_{i=-M+1}^{M} \eta_{2,i}^{el}
  \end{align}
  where the computable global error estimator $\eta_2$ and the computable local
  error estimators $\eta_{2,i}^{at}$ and $\eta_{2,i}^{el}$, associated with atoms and
  elements, respectively, are defined as
  \begin{equation} \label{EqEta2Def}
  \begin{tabular}{r@{}ll}
    $\eta_2$     & \multicolumn{2}{@{}l}{$\; := \left| \vec g^{acT} R(\vec y^{ac}) \right|
                 + \|P D (\vec y^{ac} - \vec a^a)\|_{E^a}
                   \|P D \vec g^{ac}\|_{E^a},$} \\[.5em]
    $\eta_{2,i}^{at}$ & $\; := \left| g^{ac}_i R(\vec y^{ac})_i \right|,$
                 & $\qquad i=-M+1,\ldots,M,$ \\[.5em]
    $\eta_{2,i}^{el}$ & \multicolumn{2}{@{}l}{$
                 \;:= \half \left| \big( P D (\vec y^{ac} - \vec a^a) \big)_i
                      \big( (E^a-E^{ac}) D (\vec y^{ac} - \vec a^a) \big)_i
                      \right|$} \\[.5em]
               & $\qquad + \half \left| (P D \vec g^{ac})_i
                    \big( (E^a-E^{ac}) D \vec g^{ac} \big)_i \right|,$
               & $\qquad i=-M+1,\ldots,M.$
  \end{tabular}
  \end{equation}
\end{theorem}

\section{Numerics}  \label{SecNumerics}

Now we use the two {\em a posteriori} error estimators given in
Section~\ref{SecErrorEst} to formulate an algorithm which
adaptively decides between atomistic and continuum modeling.
Finally, we present and discuss the numerical results for the periodic
array of dislocations described by the Frenkel-Kontorova model.

The error estimator $\eta_1$ gives a better estimate of the error because
$\eta_2$ involves additional inequalities. However, $\eta_2$ allows for an
atom-wise decomposition, whereas $\eta_1$ does not. This is due to the fact that
$(\eta_{low}^\pm)^2$ in the definition of $\eta_1$ is equal to the square of a
sum of atom-wise components and not the sum of the square of these components.
We can thus
let $\eta_1$ decide whether a given global tolerance $\tau_{gl}$ is already
achieved or not and use the decomposition of $\eta_2$ to decide where
the atomistic model is needed for a better approximation.  In this way, we combine
the better efficiency of $\eta_1$ with the error localization of $\eta_2$.

We start with a fully continuum model.  We then switch to the atomistic model
wherever the local error exceeds an atom-wise error tolerance $\tau_{at}$.
While decreasing $\tau_{at}$, the algorithm adaptively tags larger and larger
regions atomistic until the estimate for the goal-oriented error finally reaches $\tau_{gl}$.
The complete algorithm reads:
\begin{enumerate}
\item[(1)] Choose $\tau_{gl}$. Model all atoms as a continuum. Set $\tau_{at} \leftarrow \tau_{gl}$.
\item[(2)] Solve primal problem \eqref{EqPrimal} for $\vec y^{ac}$ and dual problem
           \eqref{EqDual} for $\vec g^{ac}$.
\item[(3)] Compute error estimator $\eta_1$ from \eqref{EqEta1Def}.
\item[(4)] If $\eta_1 \le \tau_{gl},$ then stop.
\item[(5)] Compute local error estimators $\eta_{2,i}^{at}$ and $\eta_{2,i}^{el}$
           from \eqref{EqEta2Def}.
\item[(6)] Set $\tau_{at} \leftarrow \frac{\tau_{at}}{\tau_{div}}$.
\item[(7)] Make all atoms $i$ atomistic ($\delta^a_i=1$ and $\delta^c_i=0$) for which
  \begin{align} \label{EqAdaptCrit}
    \eta_{2,i}^{tot}:=\eta_{2,i}^{at} + \half \big( \eta_{2,i-1}^{el} + \eta_{2,i}^{el} \big) \ge \tau_{at}.
  \end{align}
\item[(8)] Go to (2).
\end{enumerate}
The factor $\tau_{div}>1$ describes the rate at which the atom-wise tolerance $\tau_{at}$
is decreased during the adaptive process. We found that $\tau_{div}=10$ gives an
efficient method for this problem.

Now we come to the results for the Frenkel-Kontorova dislocation model for a periodic
chain of 1000 atoms, that is $M=500$.  The elastic constants are set to be $k_0=1$ and
$k_1=k_2=2$. To define the quantity of interest, we choose the average
displacement of atoms $11\ldots30$.  This leads to
\begin{align}
  \vec q = (q_i)_{i=-M+1,\ldots,M}, \qquad
  q_i = 1 \quad \text{for } 11 \le i \le 30, \qquad
  q_i = 0 \quad \text{otherwise.}
\end{align}
The global tolerance is chosen to be $\tau_{gl}=10^{-10}$.

\begin{table}
\begin{tabular}{|c|c|c|c|}
\hline
iteration & atomistic region & $\tau_{at}$ & $\eta_1$ \\
\hline
   1 & none                  & 1.000000e-10 & 6.860546e-03 \\
   2 & $-26 \: \ldots \: 55$ & 1.000000e-11 & 1.238016e-07 \\
   3 & $-30 \: \ldots \: 60$ & 1.000000e-12 & 2.600112e-08 \\
   4 & $-34 \: \ldots \: 66$ & 1.000000e-13 & 3.922946e-09 \\
   5 & $-38 \: \ldots \: 73$ & 1.000000e-14 & 4.104868e-10 \\
   6 & $-43 \: \ldots \: 80$ & 1.000000e-15 & 4.105166e-11 \\
\hline
\end{tabular}
\caption{\label{TabConvergence} Convergence of the algorithm for $\tau_{gl}=10^{-10}$.}
\end{table}

Table~\ref{TabConvergence} shows how the successive adaptive determination
of the atomistic-continuum modeling
proceeds. After six iterations, the atom-wise tolerance is small enough so that
$\eta_1 \le \tau_{gl}$, that is the desired accuracy has been achieved.

\begin{figure}
\includegraphics[width=0.5\textwidth]{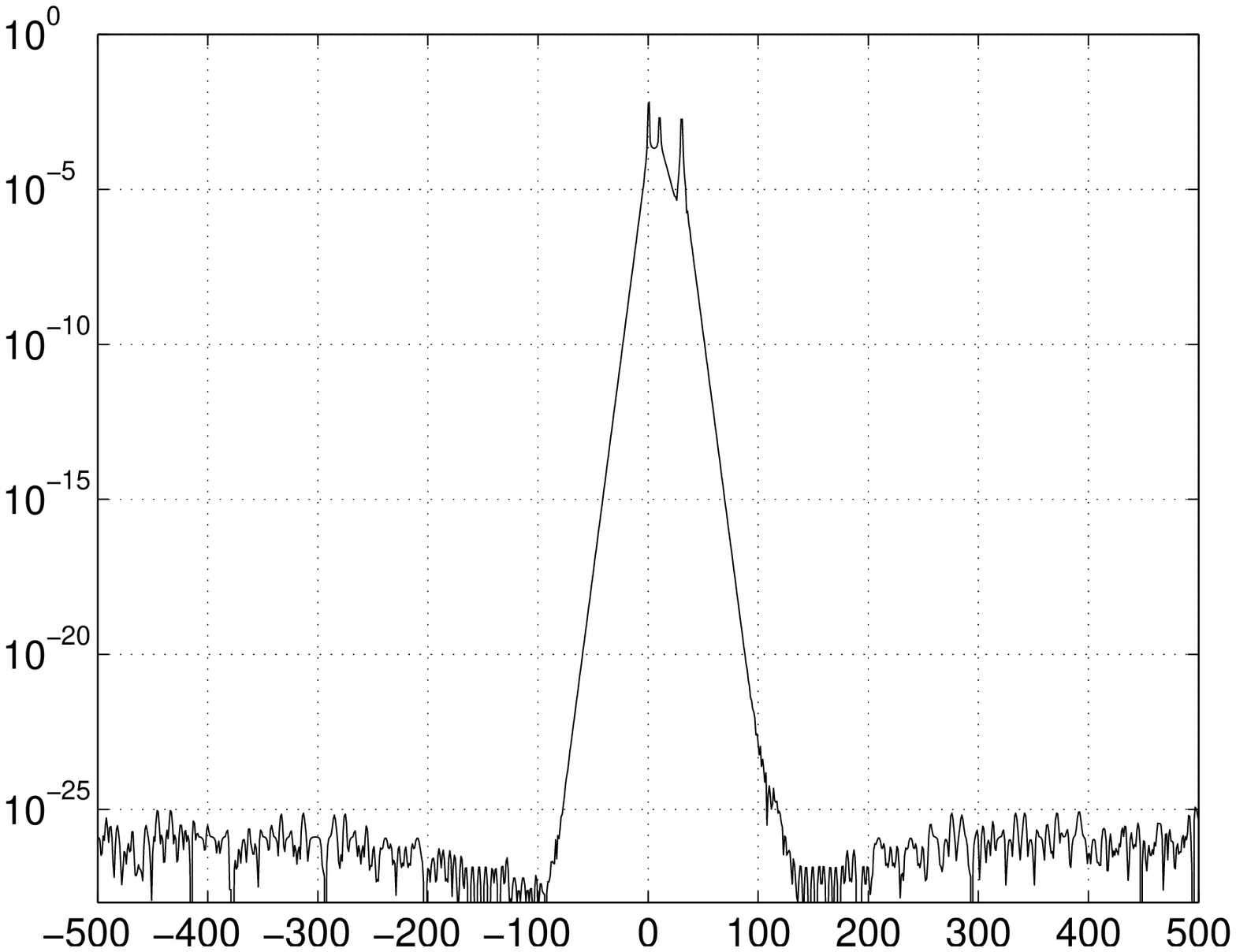}%
\includegraphics[width=0.5\textwidth]{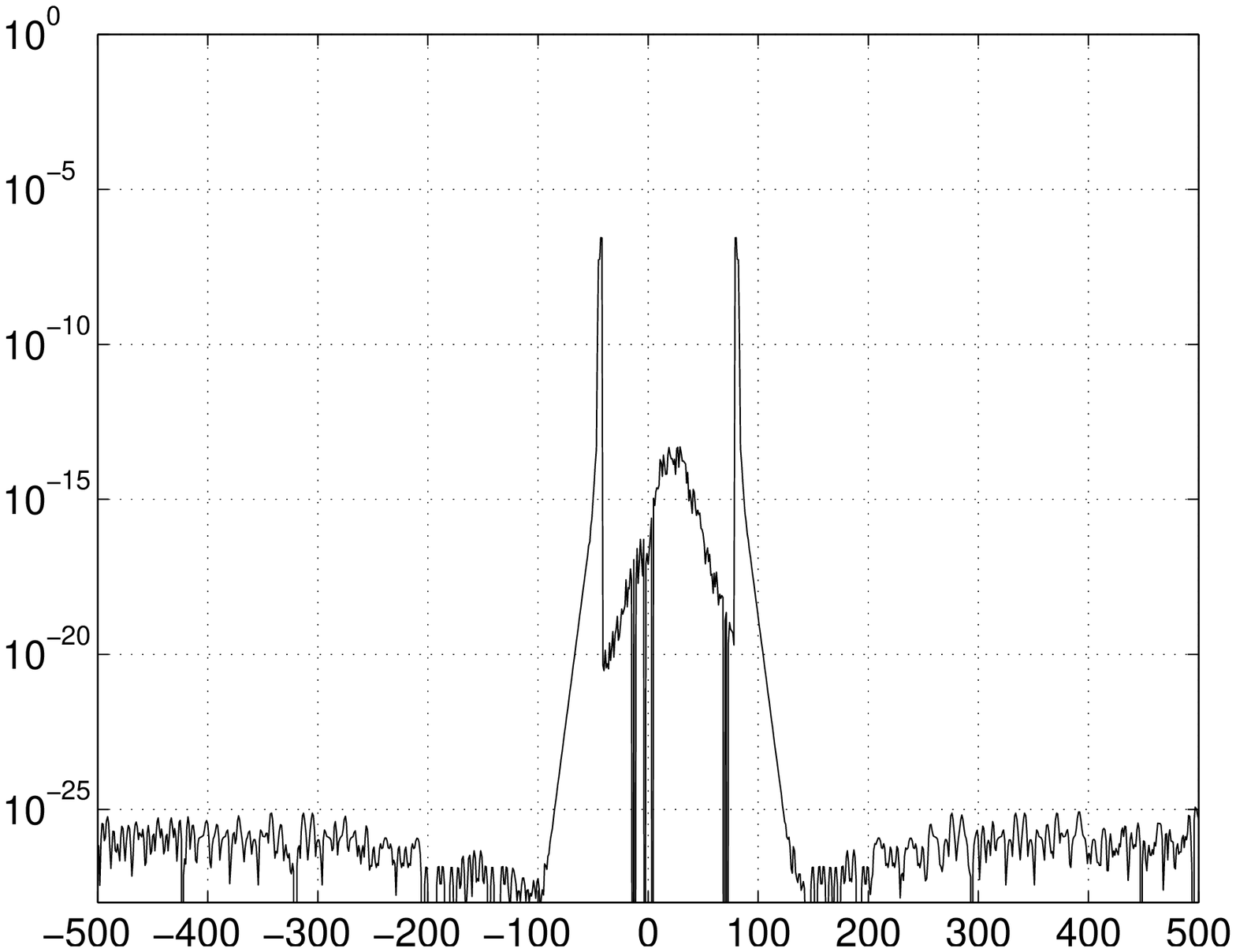}
\caption{Decomposition of the error estimator $\eta_{2,i}^{tot}$
for iteration~1 (left, fully continuum model) and for iteration~6
(right, atomistic region $-43\ldots80$).}
\label{FigAtomWise}
\end{figure}

Figure~\ref{FigAtomWise} shows the decomposition $\eta_{2,i}^{tot}$ of the error
estimator $\eta_2$ for the fully continuum model in iteration~1 of the adaption
process. One can clearly read off from the graph that the error is large near the
dislocation between atoms 0 and 1 and near atoms 11 and 30, and that it decays
exponentially away from these points. We note the slight
nonsymmetry of the atomistic-continuum modeling
due to using a goal function which averages over atoms $11\ldots30$ to the
right of the dislocation, but not to its left. The graph on the right shows the
decomposition of $\eta_2$ in the final iteration~6 with an atomistic region given by indices
$-43\dots80$. It exhibits the same nonsymmetry, but the error is considerably
smaller with peaks at the boundary between the atomistic region and the
continuum region.  In both diagrams, the fluctuations come from the limited
relative machine precision of about $10^{-16}$.

\begin{table}
\begin{tabular}{|r@{}c@{}l|c|c|c|c|c|}
\hline
\multicolumn{3}{|c|}{atomistic} &&&&& \\
\multicolumn{3}{|c|}{region}
& \raisebox{0.5em}[0mm][0mm]{$|Q(\vec y^a-\vec y^{ac})|$}
& \raisebox{0.5em}[0mm][0mm]{$\eta_1$}
& \raisebox{0.7em}[0mm][0mm]{$\displaystyle \frac{\eta_1}{|Q(\vec y^a-\vec y^{ac})|}$}
& \raisebox{0.5em}[0mm][0mm]{$\eta_2$}
& \raisebox{0.7em}[0mm][0mm]{$\displaystyle \frac{\eta_2}{|Q(\vec y^a-\vec y^{ac})|}$} \\
\hline
\multicolumn{3}{|c|}{none}    & 1.416421e-03 & 6.860545e-03 & 4.843577 & 1.231314e-02 & 8.693133 \\
  -4 & $\:\:\ldots\:\:$ &  10 & 1.863104e-03 & 6.107510e-03 & 3.278136 & 1.049800e-02 & 5.634680 \\
  -9 & $\:\:\ldots\:\:$ &  20 & 1.000572e-05 & 3.358722e-04 & 33.56803 & 6.621488e-04 & 66.17705 \\
 -14 & $\:\:\ldots\:\:$ &  30 & 1.430363e-04 & 3.187552e-04 & 2.228492 & 5.140285e-04 & 3.593694 \\
 -19 & $\:\:\ldots\:\:$ &  40 & 1.675490e-05 & 2.626711e-05 & 1.567727 & 3.691344e-05 & 2.203142 \\
 -24 & $\:\:\ldots\:\:$ &  50 & 7.361419e-07 & 1.190138e-06 & 1.616723 & 1.693910e-06 & 2.301065 \\
 -29 & $\:\:\ldots\:\:$ &  60 & 3.139276e-08 & 5.157753e-08 & 1.642975 & 7.388556e-08 & 2.353586 \\
 -34 & $\:\:\ldots\:\:$ &  70 & 1.146997e-09 & 2.001550e-09 & 1.745035 & 2.934377e-09 & 2.558312 \\
\hline
\end{tabular}
\caption{\label{TabEfficiency} Efficiency of the error estimators, $\eta_1/|Q(\vec y^a-\vec y^{ac})|$
  and $\eta_2/|Q(\vec y^a-\vec y^{ac})|$.}
\end{table}

Finally, Table~\ref{TabEfficiency} shows the efficiency of the error estimators
$\eta_1$ and $\eta_2$ for different atomistic regions. $|Q(\vec y^a-\vec
y^{ac})|$ gives the actual error which can be computed for this relatively small
problem.  In real applications, it is of course not available.  One can clearly
see that $\eta_1$ gives a better estimate than $\eta_2$, which numerically
confirms our conjecture that $\eta_1$ is a better estimator than $\eta_2$.  An
unusually high value for the efficiency occurs when the atomistic-continuum
boundary sweeps through the region where the quantity of interest is measured.
After this, the efficiencies converge to decent values around 1.7 and 2.5 for
$\eta_1$ and $\eta_2$, respectively. We note that for clamped boundary
conditions and a symmetric quantity of interest, better efficiencies of 1.4 and
2, respectively, have been obtained \cite{ArndtLuskin:2007a}.

\bibliographystyle{hsiam}
\bibliography{../Literature/marcel}

\end{document}